\documentclass[11pt]{amsart}
\usepackage{bbm,pstricks,amsmath,amssymb,latexsym,fullpage}
\usepackage{epsf,psfrag,epsfig}
\usepackage[english]{babel}
\usepackage[latin1]{inputenc}

\newtheorem{thm}{Theorem}

\newtheorem{lemma}[thm]{Lemma}

\newcommand{\od}{\!\odot\!}
\newcommand{\op}{\!\oplus\!}

\newcommand{\dis}{\displaystyle }

\renewcommand{\Pr}{\mathbb{P}}

\title{Short proof of Rayleigh's Theorem and extensions}

\author{Olivier Bernardi}
\thanks{I acknowledge partial support from ANR A3 and European project ExploreMaps.}

\begin{document}

\begin{abstract}
Consider a walk in the plane made of $n$ unit steps, with directions chosen independently and uniformly at random at each step. Rayleigh's theorem asserts that the probability for such a walk to end at a distance less than 1 from its starting point is $1/(n+1)$. We give an elementary proof of this result. We also prove the following generalization valid for any probability distribution $\mu$ on the positive real numbers: if two walkers start at the same point and make respectively $m$ and $n$ independent steps with uniformly random directions and with lengths chosen according to $\mu$, then the probability that the first walker ends farther than the second is~$m/(m+n)$. 
\end{abstract}

\maketitle

\begin{figure}[h!]
\centerline{\includegraphics[scale=.7]{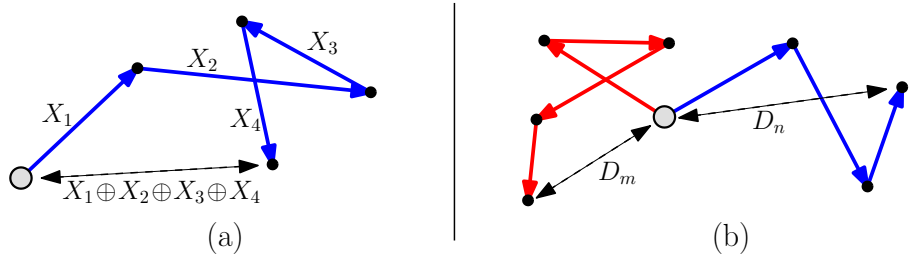}}
\caption{(a) The distance $X_1\oplus X_2\oplus X_3\oplus X_4$ achieved after four steps. 
(b)~Comparing the distances $D_m=m\od X$ and $D_n=n\od X$.}
\label{fig:walk}
\end{figure}

We consider random walks in the Euclidean plane. Given some real positive random variables $X_1,X_2,\ldots,X_n$, we consider a random walk starting at the origin of the plane and made of $n$ steps of respective length $X_1,X_2,\ldots,X_n$, with the direction of each step chosen independently and uniformly at random. We denote by $X_1\oplus X_2\oplus \cdots \oplus X_n$ the random variable corresponding to the distance between the origin and the end of the walk. This definition is illustrated in Figure~\ref{fig:walk}(a).

For a non-negative real random variable $X$, we denote by $n\od X$ the random variable $X_1\oplus \cdots \oplus X_n$, where $X_1,\ldots, X_n$ are independent copies of $X$. Hence $n\od X$ represent the final distance from the origin after taking $n$ independent steps of lengths distributed like $X$ and directions chosen uniformly at random.
Rayleigh's theorem asserts that if $X=1$, that is, each step has unit length, then for all $n>1$,
$$\Pr(n\od X <\, 1)=\frac{1}{n+1}.$$
This theorem was first derived from Rayleigh's investigation of ``random flights'' in connection with Bessel functions (see \cite{Watson:Bessel-functions}) and appears as an exercise in \cite[p.104]{Spitzer:random-walks}\footnote{The exercise calls for developing the requisite Fourier analysis for spherically symmetric functions in order to obtain an identity involving Bessel functions.}. A simpler proof was given by Kenyon and Winkler as a corollary of their result on branched polymers~\cite{Kenyon:branched-polymers}. The goal of this note is to give an elementary proof of the following generalization of Rayleigh's theorem. 



\begin{thm}\label{thm:extension-Rayleigh}
Let $X$ be a real random variable taking positive values, and let $m,n$ be non-negative integers such that $m+n>2$. If $D_m$ and $D_n$ are independent random variables distributed respectively like $m\!\odot\! X$ and $n\!\odot\! X$, then 
$$\Pr(D_m>D_n)=\frac{m}{m+n} .$$
In words, if two random walkers start at the origin and take respectively $m$ and $n$ independent steps with uniformly random directions and with lengths chosen according to the distribution of~$X$, then the probability that the first walker ends farther from the origin than the second walker is $m/(m+n)$.
\end{thm}

Theorem~\ref{thm:extension-Rayleigh} is illustrated in Figure~\ref{fig:walk}(b). Clearly, this extends Rayleigh's theorem which corresponds to the case $m=1$ and $X=1$. Our proof of Theorem~\ref{thm:extension-Rayleigh} starts with a lemma based on the fact that the angles of a triangle sum to~$\pi$.

\begin{lemma}\label{lem:triangle}
For any random variables $A,B,C$ taking real positive values, 
\begin{equation}\label{eq:triangle}
\Pr(A>B\oplus C)+\Pr(B>A\oplus C)+\Pr(C>A\oplus B)=1.
\end{equation}
\end{lemma}

\begin{figure}[h!]
\centerline{\includegraphics[scale=.7]{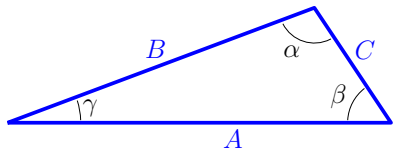}}
\caption{The triangle $T$ with side lengths $A,B,C$.}
\label{fig:triangle}
\end{figure}

\begin{proof}
By conditioning on the values of the random variables $A,B,C$, it is sufficient to prove \eqref{eq:triangle} in the case where $A,B,C$ are non-random positive constants, and the randomness only resides in the directions of the steps. 
Now we consider two cases. First suppose that one of the lengths $A,B,C$ is greater than the sum of the two others. 
In this case, one of the probabilities appearing in \eqref{eq:triangle} is 1 and the others are 0, hence the identity holds. 
Now suppose that none of the lengths $A,B,C$ is greater than the sum of the two others. In this case, there exists a triangle $T$ with side lengths $A,B,C$. The triangle $T$ is shown in Figure \ref{fig:triangle}. 
The probability $\Pr(A>B\op C)$ is equal to $\alpha/\pi$, where $\alpha$ is the angle between the sides of length $B$ and $C$ in the triangle $T$ (because $A>B\oplus C$ if and only if the angle between the step of length $B$ and the step of length $C$ is less than $\alpha$ in absolute value). Summing this relation for the three probabilities appearing in \eqref{eq:triangle} gives
$$\Pr(A>B\oplus C)+\Pr(B>A\oplus C)+\Pr(C>A\oplus B)=\frac{\alpha+\beta+\gamma}{\pi}=1.$$
where $\alpha,\beta,\gamma$ are the angles appearing in Figure \ref{fig:triangle}.
\end{proof}

We now complete the proof of Theorem~\ref{thm:extension-Rayleigh}. 
Let $s=m+n$ and let $D_0,D_1,\ldots,D_s$ be independent random variables distributed respectively like $0\od X,1\od X,\ldots, s\od X$. We denote $p_i=\Pr(D_i>D_{s-i})$ and want to prove $\dis p_m=m/s$. Let $i,j,k$ be positive integers summing to $s$. Applying Lemma~\ref{lem:triangle} to $A=D_i$, $B=D_j$, $C=D_k$ gives $p_i+p_j+p_k=1$. 
Moreover, $p_k=1-p_{s-k}$ since $\Pr(D_k=D_{s-k})=0$ (recall that $s>2$). 
Thus 
$$p_i+p_j=p_{i+j},$$
for all $i,j>0$ such that $i+j\leq n$. By induction, this implies 
$i\, p_1=p_i$ for all $i\in \{1,\ldots,n\}$. In particular $p_1=p_s/s=1/s$, and $p_m=m\,p_1=m/s$. This concludes the proof of Theorem~\ref{thm:extension-Rayleigh}.\\

\smallskip

\noindent \textbf{Acknowledgments:} I thank Peter Winkler for extremely stimulating discussions.


\end{document}